\documentclass[showpacs,showkeys,amssymb,aps,notitlepage]{revtex4-1}

\usepackage{dcolumn}
\usepackage{bm}

\usepackage{amsmath}
\usepackage{etoolbox} 

\usepackage[pdfauthor={Richard J. Mathar}, colorlinks=true,
 pdfkeywords={Java BigDecimal Functions},bookmarksopen=true,pdfpagelayout=OneColumn]{hyperref}

\usepackage{orcidlink}
\begin{document}

\title[Java BigDecimalMath]{A Java Math.BigDecimal Implementation of Core Mathematical Functions}

\author{Richard J. Mathar \orcidlink{0000-0001-6017-6540}}
\homepage{https://mathar.www3.mpia.de/}
\affiliation{Max Planck Institute for Astronomy, K\"onigstuhl 17, 69117 Heidelberg, Germany}

\date{\today}

\begin{abstract}
The mathematical functions $\log(x)$, $\exp(x)$, $\sqrt[n]x$, $\sin(x)$, $\cos(x)$,
$\tan(x)$, $\arcsin(x)$, $\arctan(x)$, $x^y$, $\sinh(x)$, $\cosh(x)$, $\tanh(x)$ and $\Gamma(x)$ have been
implemented for arguments
$x$ in the real domain in a native Java library
on top of the multi-precision \texttt{BigDecimal} representation of floating point numbers. This supports
scientific applications where more than the double precision accuracy of the library
of the Standard Edition is desired. The full source code is made available under the LGPL v3.0.
\end{abstract}

\maketitle

\section{Overview}
\subsection{Aim}
Whereas many Java applications can use the Java Native Interface
to bind to (C-based) multi-precision programs on a \emph{host} platform
when a higher precision
than the 64-bit standard is needed \cite{FousseTOMS33}, others may observe that there is
only rudimentary
support on the \emph{native} platform if standard mathematical functions are needed.

The aim of this script is to provide a base implementation
of the core trigonometric and algebraic functions \cite{SmithTOMS17}
on top of the native
\texttt{BigDecimal} class with infinite-precision capability. Demand
originates from scientific and perhaps engineering computations, where accumulation
of rounding errors (loss of digits) might pose a problem.

\subsection{Design Choices}

A characteristic feature of the implementation suggested here is that
the floating point variables which are arguments to mathematical functions
define by their number of digits
which precision is
achieved in the result. The estimate of the accuracy of the result
is generally derived from the first order Taylor approximation of the function in question
based on the accuracy of the input variable. The number of digits of the
values returned will be larger than the number of digits of the variable where
the function is flat---for example the $\arctan(x)$ where $x\gg 1$---, smaller where it is steep.
This is a deliberate
difference to most computer algebra systems. It provides some semi-automate
detection of loss-of-precision through cancellation of digits, and it 
reduces some burden to the application programmer to decide on a
\texttt{MathContext} interface prior to each individual call.

This is backed by some \texttt{xxxRound} functions where \texttt{xxx}
are the fundamental \texttt{add}, \texttt{subtract}, \texttt{multiply},
\texttt{divide} operations, etc.,
which internally calculate estimators of the precision of the result based on
the precisions of their arguments. The mathematical functions are typically
some power series expansions, and they make heavy use of these to keep
the error accumulation of the individual terms under control, that is, to
chop off digits early to keep execution times short where intermediate results
are known to be dominated by noise in the parameters.

As a side effect, the number of digits returned may be even smaller
than the characteristic 6 digits of a single-precision calculation.
In addition, the results depend on the number of trailing zeros of the inputs.
(A function \texttt{scalePrec} is provided to boost the apparent accuracy
of numbers by appending zeros.)

\subsection{Known Limitations}
As presented, the implementation is known to have jitters of 1 or 2 
in the least significant digits in some values returned.

The algorithms have been chosen for reliability and simplicity,
and may be
slower than alternatives which have not been investigated.

Classes of important special functions (Bessel functions, Elliptic Integrals,\ldots)
and complex arithmetic are absent.

\section{Implementation Strategies}

\subsection{Constants}

The heavy-duty constants $\pi$, $e$, $\ln 2$ and $\gamma$
are tabulated to high precision which presumably suffices
for most purposes in engineering and sciences. The backup implementations
for applications in some areas of mathematics are:
\begin{itemize} 
\item
$\pi$ is evaluated by
Broadhurst's
equation (18) \cite{BroadhurstArxiv98}.
\item
$\log 2$ is evaluated by
Broadhurst's
equation (21) \cite{BroadhurstArxiv98} finalized by pulling
a square root.
\item
$\gamma$ is generated  by the series \cite[(3.9)]{DilcherAM48}
\begin{equation}
\gamma = 1-\log\frac{3}{2}-\sum_{n\ge 1}\frac{\zeta(2n+1)-1}{4^n(2n+1)}
.
\end{equation}
\item
$e$ is forwarded to the generic evaluation of $\exp(1)$, Section \ref{sec.exp}.
\end{itemize} 

\subsection{Roots}\label{sec.root}

The roots $y=\sqrt[n]x$ for positive integer $n$, including the special case
of square roots $n=2$, are computed by iterative updates with
the first order Newton method \cite[(3.96)]{AS}
\begin{equation}
y\to y-\frac{1}{n}\left( y-\frac{x}{y^{n-1}} \right)
.
\end{equation}
The initial estimates of $y$ are
set in double precision by a call to the \texttt{Math.pow}.

The \texttt{hypot} function computes
\begin{equation}
z=\sqrt{x^2+y^2}
\end{equation}
from two arguments $x$ and $y$. A derived case has been implemented
taking an integer value $x$, because this implies that the precision of the
result $z$ is determined from the precision in $y$ alone.

\subsection{Exponential}\label{sec.exp}
If $x$ is close to zero, the standard Taylor series \cite[(4.2.1)]{AS}
\begin{equation}
\exp(x) = \sum_{k\ge 0} \frac{x^k}{k!}
\label{eq.expT}
\end{equation}
is employed. For larger $x$, $x$ is scaled by powers of 10
\begin{equation}
e^x=\left( e^{10^{-t}x} \right)^{10^t}
\end{equation}
such that the value in parenthesis can be evaluated by recourse to (\ref{eq.expT}).
Scaling by powers of 10 is a cheap operation in the \texttt{BigDecimal}
library because it only involves a diminuation of the scale.
Powers with integer exponents are also more efficient than one might na\"{\i}vely expect. The
10th power needs 4 multiplications, for example; see \cite{McCarthyMathComp46}
and sequence A003313 in the OEIS \cite{sloane}.

The general powers are forwarded to a mixed call of the \texttt{log}
and \texttt{exp} functions,
\begin{equation}
x^y = \exp( y\log x).
\end{equation}

\subsection{Logarithm}
For arguments close to 1, the standard Taylor expansion
\cite[(4.1.24)]{AS}
\begin{equation}
\log(1+x)=\sum_{k\ge 1}(-1)^{j+1}\frac{x^k}{k}
\end{equation}
is used. For larger $x$, adaptation to that range is achieved by
scaling with some integer $r$
\begin{equation}
\log x = r\log \sqrt[r]x
\end{equation}
with an auxiliary call to the \texttt{root} function of Section \ref{sec.root}.
The variable $r$ is obtained by a call to the \texttt{Math.log} of the
native  library.

For some integer arguments, dedicated routines are implemented assuming
that $\ln 2$ is instantly available,
\begin{eqnarray}
12 \ln 3 &=& 19 \ln 2 +\sum_{k\ge 1}\frac{(-1)^{k+1}}{k}\left(\frac{7153}{524288}\right)^k, \\
6 \ln 5 &=& 14 \ln 2 -\sum_{k\ge 1}\frac{1}{k}\left(\frac{759}{16384}\right)^k, \\
\ln 7 &=& 3 \ln 2 -\sum_{k\ge 1}\frac{1}{k8^k}.
\end{eqnarray}
These have practically no speed advantage compared to the alternative of
adding zeros and handling them with the generic procedure described above.

\subsection{Trigonometric}
The arguments of $\sin$, $\cos$ and $\tan$ are reduced to the fundamental
domain modulo $2\pi$ or modulo $\pi$, then folded with standard shifting
equations, Table 4.3.44 in the Handbook \cite{AS} into regions where
the basic Taylor series converge well. These are in particular
\begin{eqnarray}
\sin x &=& \sum_{k\ge 0} (-1)^k\frac{x^{2k+1}}{(2k+1)!}, \\
\cos x &=& \sum_{k\ge 0} (-1)^k\frac{x^{2k}}{(2k)!}
\end{eqnarray}
for $x< \pi/4$, and
\begin{equation}
\tan x = \sum_{k\ge 1} (-1)^{k+1}\frac{4^k(4^k-1)}{(2k)!}B_{2k} x^{2k-1}
\label{eq.tanT}
\end{equation}
for $x<0.8$ \cite[(4.3.67)]{AS}. The $\tan x$ is forwarded to a similar
expansion of $\cot x$ \cite[(4.3.70)]{AS} if $x>0.8$.

\subsection{Inverse Trigonometric}
The $\arcsin$ is implemented as \cite[(4.4.40)]{AS}
\begin{equation}
\arcsin x = \sum_{k\ge 0} \frac{(2k-1)!!}{(2k)!! (2k+1)}x^{2k+1}
\end{equation}
where $x<0.7$, and as the complementary \cite[(4.4.41)]{AS}
where $0.7<x<1$. The $\arctan$ is implemented by the
standard Taylor expansions
\begin{eqnarray}
\arctan x = \sum_{k \ge 0} (-1)^k \frac{x^{2k+1}}{2k+1},\quad x<0.7\\
\arctan x = \frac{\pi}{2}-\sum_{k \ge 0} (-1)^k \frac{1}{(2k+1)x^{2k+1}},\quad x>3.
\end{eqnarray}
The intermediate cases are mapped to the region $x<0.7$ by reverse application
of \cite[(4.4.34)]{AS}
\begin{equation}
2\arctan x = \arctan\frac{2x}{1-x^2}
\end{equation}
at the cost of one additional square root.

\subsection{Hyperbolic}

The hyperbolic functions $\sinh$ and $\cosh$ are evaluated by
their power series \cite[(4.5.62),(4.5.63)]{AS} if the argument is
close to zero, otherwise transformed by the multi-angle formulas. The
$\tanh$ is implemented as
\begin{equation}
\tanh x = \frac{1-\exp(-2x)}{1+\exp(-2x)}.
\end{equation}
The inverse hyperbolic functions are mapped to their logarithmic representations.

\subsection{Gamma Function}
The $\Gamma$ function is reduced to the region near $x=1$ by its
functional equation
\begin{equation}
x\Gamma(x)=\Gamma(x+1),
\end{equation}
and then expanded with \cite[(6.1.33)]{AS}
\begin{equation}
\ln \Gamma(1+x) = -\ln(1+x)+x(1-\gamma) +\sum_{k\ge 2}(-1)^k\frac{\zeta(k)-1}{k}x^k
.
\end{equation}
This bypasses difficulties of regulating the errors in the Stirling formula
\cite{SpiraMathComp25,WrenchMCom22,GordonJACM7}, but needs a rather costly evaluation of the $\zeta$ function.
For even indices we implement \cite[(23.2.16)]{AS}
\begin{equation}
\zeta(2n)= \frac{(2\pi)^{2n}}{2(2n)!}|B_{2n}|,
\label{eq.zetaB}
\end{equation}
for indices $3$ or $5$ the Broadhurst expansions \cite{BroadhurstArxiv98},
and for odd arguments $\ge 7$ \cite{CohenExpMath1,VepstasRJ27}
\begin{eqnarray}
\zeta(n) = \frac{(2\pi)^n}{n-1}
\sum_{k=0}^{(n+1)/2}
(-1)^k(1-2k)\frac{B_{2k}B_{n+1-2k}}{(2k)!(n+1-2k)!}
\nonumber \\
-2\sum_{k\ge 1}
\frac{1}{k^n(e^{2\pi k}-1)}
\left(1+\frac{2\pi k\epsilon_n}{1-e^{-2\pi k}}\right),
\end{eqnarray}
where
\begin{equation}
\epsilon_n = \left\{
\begin{array}{ll}
0, & n\equiv 3 \pmod 4,\\
2/(n-1), & n\equiv 1 \pmod 4.
\end{array}
\right.
\end{equation}

\section{Integer Classes}

The aim to delay rounding of rational numbers leads to the auxiliary
implementation of a \texttt{Rational} data type which consists of a
signed numerator and an unsigned denominator, both of the \texttt{BigInteger} type.
The basic operations of multiplication, division, addition and raising
to an integer power are all exact in that class, and also some 
integer roots if numerator and denominator are perfect powers.

The class \texttt{Bernoulli} creates a special instance of these
rational numbers, the Bernoulli numbers which are helpful in
(\ref{eq.tanT}) and (\ref{eq.zetaB}).
From a short initial table at small indices \cite[(Tab 23.2)]{AS},
values at larger indices are generated by a double
sum \cite[(1)]{GouldAMM79}:
\begin{equation}
B_n = \sum_{k=0}^n \frac{1}{k}\sum_{j=0}^k (-1)^j j^n \binom{k}{j}
.
\end{equation}

This is augmented by a very rudimentary \texttt{Prime}
class which grows dynamically, and a class \texttt{Ifactor}
which represents a positive integer and its prime number factorization.
The multiplicative sums-of-divisors function $\sigma_k()$ and $\varphi()$
of number theory are derived from
such an intermediate prime number decomposition. The set of divisors
of an integer is created as well from there by a multinomial scan
of the exponents.

Numbers of the form $r\sqrt d$, where $r$ is a signed rational number
and $d$ a non-negative rational number of the \texttt{Rational} class,
are represented by the \texttt{BigSurd}
class. Sums of these are represented by \texttt{BigSurdVec}, for which
addition, subtraction and exponentiation with integer exponents can
be represented exactly.

\subsection{Application: Wigner $3n$-$j$ symbols}
As an independent test of of other programs that evaluate Wigner 3n-$j$
symbols \cite{MatharArxiv1102a}, footed on the exact representation
of square roots presented above, the class \texttt{Wigner3j} allows
computation of 3n-$j$ symbols for unlimited $n$.
The application interface has been modeled according to an earlier
proposal \cite{BarCPC50}. One line is the integer $3n$. Two lines of integers are 1-based
indices into a list of $j$-values, implicitly bundled in triads such
that for each of the $2n$ factors of the underlying cubic graph the three
contributing $j$-values are listed in the order of appearance in the
upper row of their $3jm$ values. These indices are negated if
the associated $m$-value in the $3jm$ symbol appears with a negative
sign.
All further lines of the input contain lists of $j$ values by actually
providing the positive, integer-valued $2j+1$.

The implementation does not use any reduction techniques. It performs
the summation over all $m$-values of all angular momenta that have been
defined, and is slow for that reason. (There is some look-ahead for
small cycles in the associated cubic graph to take advantage of the
selection rule of $3jm$-values.)

The class \texttt{Wigner3jGUI} is an online calculator for these.
The connectivity schemes for all symbols from $6j$ up of $15j$
(in a serialized order reading the entries by rows in the braces
of the standard notation) have been initialized and are selectable with a button.

\subsection{Application: Fitting to P-recursive sequences}
\subsubsection{(Finite) Integer Sequence}
The class \texttt{GuessPRec} has constructors
that store a (finite) list/sequence of integer numbers $a(n)$ defined 
over a set of integers $n$---the minimum index $n$ will
be called the \textit{offset} following an OEIS convention \cite{sloane}.
The class has a function \texttt{guess} that
tries to find a list of polynomials $p_i(n)$ solving the equation
\begin{equation}
\sum_{i=0}^o p_i(n)a(n-i) + p_{o+1}(n) =0
\label{eq.P}
\end{equation}
with a small degree of the polynomials $p_i$ and a small order $o$ of the recurrence.
Polynomials of degree 1 (constants, C-finite fits) are admitted, and we also admit solutions
where the polynomial $p_{o+1}(n)$ is zero, i.e., where the recurrence is \textit{homogeneous},
or where some of the intermediate $p_{n-1},\ldots, p_{o-1}$ are zero.

For some types of sequences it is known a priori that P-finite recurrences
exist \cite{KauersTetrah,MatharArxiv2109}.
For sequences that enumerate objects that can be computed by the transfer matrix method,
the generating functions are rational polynomials which lead to C-finite recurrences
\cite[\S V.6]{Flajolet}\cite[\S 3.3.3]{Flajolet4103}.

\subsubsection{Linear Algebra}
The search space of the algorithm considers a fixed number $o$ of terms
with a fixed degree $d$ of the polynomials and adaptive integer coefficients $c$:
\begin{equation}
p_i(n) = \sum_{l=0}^d c_{i,l} n^l.
\end{equation}
With that ansatz there are $N=(o+2)(d+1)$ free/unknown coefficients $c$, and 
the long writeup of \eqref{eq.P} 
is
\begin{equation}
[c_{0,0}+c_{0,1}n+c_{0,2}n^2+\cdots]a(n)
+
[c_{1,0}+c_{1,1}n+c_{1,2}n^2+\cdots]a(n-1)
+\cdots
c_{o+1,0}+c_{o+1,1}n+c_{o+1,2}n^2+\cdots =0,
\label{eq.Plong}
\end{equation}
which involves $o+1$ coefficients $a(n)$ down to $a(n-o)$.
We require that this equation is true for some starting index $n$ down to an 
index $n-N+1$,
and implement these requirements by a set of $N=(o+2)(d+1)$ linear equations
\begin{eqnarray}
\left(
\begin{array}{cccccccccc}
a(n) & na(n) & n^2 a(n) & \hdots & a(n-1) & na(n-1) & \cdots 1 & n & \cdots & n^d \\
a(n-1) & (n-1)a(n-1) & (n-1)^2 a(n-1) & \hdots & a(n-2) & (n-1)a(n-2) & \cdots 1 & n-1& \cdots & (n-1)^d \\
\cdots \\
\end{array}
\right)
\nonumber
\\
\cdot 
\left(
\begin{array}{c}
c_{0,0}\\
c_{0,1}\\
\cdots\\
c_{1,0}\\
c_{1,1}\\
\cdots\\
c_{o+1,0}\\
c_{o+1,1}\\
\cdots\\
c_{o+1,d}\\
\end{array}
\right)
=
\left(
\begin{array}{c}
0\\
0\\
0\\
\cdots\\
0
\end{array}
\right).
\end{eqnarray}
In the coefficient matrix on the left hand side, the $n$ in a row is the $n$ of the 
previous row decreased by one.
[Due to that shift of the indices by one in this implementation, the set of the initial/known $a$
must not have holes of undefined $a$ at intermediate indices.]

Because we start from the homogeneous equation \eqref{eq.Plong}, zero on the right hand side,
the matrix is rank-deficient and
the set of the $c_{i,l}$
is defined only up to a common constant---which is conveniently defined to
minimize the least common multiple of these constants (i.e., to divide them such
that the greatest common divisor of the $c$ is one) and to flip all signs
that for example the non-vanishing $c_{0,d}$ of the largest $d$ is positive.

\subsubsection{Outer Loop}
Typically neither the polynomial degrees $d$ nor the number of terms in the recurrence
are known beforehand. So the function \texttt{guess} starts from $N=1$ and increases $N$
until the degrees of freedom are exhausted, i.e., until $N$ is larger than the number
of known integers $a$ in the sequence. For each $N$ we loop over two divisors $o+2$ and $d+1$
and test whether the linear equations of the aforementioned type emit
vectors $c$ which are constant for a maximum range of starting indices $n$.

After compilation, the algorithm is started with the \texttt{main} class of \texttt{GuessPRec} 
in the jar-file as

\texttt{java -cp org-nevec-rjm-5.315.jar org.nevec.rjm.GuessPRec} \textit{[-v]} -f \textit{filename}

where \textit{filename} is an ASCII file which defines a finite sequence of $n$ and $a(n)$ in
the following format (a weaker variant of the b-files of the OEIS \cite{sloane}):
\begin{itemize}
\item
lines that are empty or start with \# are ignored/skipped (comments)
\item
lines that start with optional white space (blanks, \ldots), an integer number $n$,
more white space, and another integer number $a(n)$ and optionally other contents
are indices $n$ and their sequence elements $a(n)$, where the first column, $n$, must
start at the offset and increase by 1 for each followup non-comment line.
\end{itemize}
The program prints lines of the format \eqref{eq.Plong} to standard output,
scanning blocks of consecutive $a(n)$ with shifting regions. The more of
these lines are the same, the more trustworthy is the conjecture that
these recurrences are valid.

\section{Summary} 
The most frequently used mathematical function
with arguments and return values of the multi-precision \texttt{BigDecimal}
type are presented in a Java library. Control over the variable
requirements in precision is basically achieved by recourse to simple algorithms
that allow semi-analytic estimations of the propagation of errors.

\appendix
\section{Compilation}
The source code is compiled into a \texttt{jar}-file
with \texttt{maven}---a package in all standard Linux systems---either 
directly via the \texttt{Makefile},
\begin{verbatim}
cd anc
make
\end{verbatim}
or better using the \texttt{autools} as
\begin{verbatim}
cd anc
autoreconf -i -f -s
./configure
make
make doc # optional: html of the API in the doc directory
\end{verbatim}
If \texttt{maven} is absent, the grassroot compilation is
\begin{verbatim}
cd anc/src/main/java
find org/nevec/rjm/*.java -print | xargs javac
jar -c -f org-nevec-rjm-5.315.jar org/nevec/rjm/*.class
mv *.jar ../../..
\end{verbatim}

\bibliographystyle{apsrmp4-1}
\bibliography{all}

\end{document}